\documentclass[12pt,a4paper]{article}

\usepackage{theorem}
\usepackage{enumerate}
\usepackage{amsmath} 
\usepackage{latexsym}
\usepackage{amssymb}
\usepackage{amsfonts}
\usepackage{graphicx}
\usepackage[varg]{txfonts}
\usepackage{here}

\theorembodyfont{\normalfont\slshape}
\newtheorem{thm}{Theorem}

\newtheorem{Thm}{Theorem}

\newtheorem{lem}{Lemma}
\newtheorem{Lem}{Lemma}

\newtheorem{claim}{Claim}

\theoremstyle{remark}

{\theorembodyfont{\upshape}
}

\newcommand{\proof}{\medbreak\noindent\textit{Proof.}\quad}

\newcommand{\qed}{{$\quad\square$\vspace{3.6mm}}}

\newcommand{\ora}{\overrightarrow}
\newcommand{\ola}{\overleftarrow}
\newcommand{\ol}{\overline}

\renewcommand{\labelenumi}{\upshape{(\theenumi)}}

\bfseries\normalfont

\begin{document}

\title{Partitioning a graph into cycles \\with a specified number of chords}

\author{Shuya Chiba\thanks{Applied Mathematics, 
Faculty of Advanced Science and Technology, 
Kumamoto University, 
2-39-1 Kurokami, Kumamoto 860-8555, Japan 
(E-mail address: \texttt{schiba@kumamoto-u.ac.jp}). 
This work was supported by JSPS KAKENHI Grant Number 17K05347.}
\and 
Suyun Jiang\thanks{School of Mathematics, 
Shandong University, 
Jinan 250100, China; 
Institute for Interdisciplinary Research, Jianghan University, Wuhan, China 
(E-mail address: \texttt{jiang.suyun@163.com}). }
\and 
Jin Yan\thanks{School of Mathematics, 
Shandong University, 
Jinan 250100, China 
(E-mail address: \texttt{yanj@sdu.edu.cn}). 
This work was supported by National Natural Science Foundation of China (No.11671232).}
}

\date{}
\maketitle

\vspace{-12pt}
\begin{abstract}
For a graph $G$, 
let $\sigma_{2}(G)$ be the minimum degree sum 
of two non-adjacent vertices in $G$. 
A chord of a cycle in a graph $G$ 
is an edge of $G$ joining two non-consecutive vertices of the cycle. 
In this paper, we prove the following result, 
which is an extension of a result of Brandt et al.~(J. Graph Theory 24 (1997) 165--173) for large graphs: 
For positive integers $k$ and $c$, 
there exists an integer $f(k, c)$ such that, 
if $G$ is a graph of order $n \ge f(k, c)$ 
and $\sigma_{2}(G) \ge n$, 
then 
$G$ can be partitioned into $k$ vertex-disjoint cycles, 
each of which has at least $c$ chords.

\medskip
\noindent
\textit{Keywords}: Partitions; Cycles; Chorded cycles; Degree sum conditions\\
\noindent
\textit{AMS Subject Classification}: 05C70, 05C45, 05C38
\end{abstract}

\section{Introduction}
\label{Introduction}

In this paper, we consider finite simple graphs, 
which have neither loops nor multiple edges.
For terminology and notation not defined in this paper, we refer the readers to \cite{BMBook2008}. 
Let $G$ be a graph. 
For a vertex $v$ of $G$, 
we denote by $d_{G}(v)$ and $N_{G}(v)$ the degree and the neighborhood of $v$ in $G$. 
Let $\delta(G)$ be the minimum degree of $G$ 
and let 
$\sigma_{2}(G)$ 
be the minimum degree sum of two non-adjacent vertices in $G$, 
i.e., 
if $G$ is non-complete, 
then $\sigma_{2}(G) = \min \big\{d_{G}(u) + d_{G}(v) : u, v \in V(G), u \neq v, uv \notin E(G) \big\}$; 
otherwise, 
let $\sigma_{2}(G) = +\infty$. 
If the graph $G$ is clear from the context, 
we often omit the graph parameter $G$ in the graph invariant. 
We denote by $K_{t}$ the complete graph of order~$t$. 
In this paper, 
``partition'' and ``disjoint'' always mean 
``vertex-partition'' 
and 
 ``vertex-disjoint'', respectively.

\medskip
A graph is \textit{hamiltonian} if 
it has a \textit{Hamilton cycle}, i.e., a cycle containing all the vertices of the graph. 
It is well known that 
determining whether a given graph is hamiltonian or not, is NP-complete. 
Therefore, 
it is natural to study sufficient conditions for hamiltonicity of graphs. 
In particular, 
since the approval of the following two theorems, 
various studies have considered degree conditions.

\begin{Thm}[Dirac \cite{Dirac1952}]
\label{thm:Dirac1952}
Let $G$ be a graph of order $n \ge 3$.
If $\delta \ge \frac{n}{2}$, 
then $G$ is hamiltonian.
\end{Thm}

\begin{Thm}[Ore \cite{Ore1960}]
\label{thm:Ore1960}
Let $G$ be a graph of order $n \ge 3$.
If $\sigma_2 \ge n$, 
then $G$ is hamiltonian.
\end{Thm}

In 1997, Brandt et al.~generalized the above theorems by showing that 
the Ore condition, i.e., the $\sigma_{2}$ condition in Theorem~\ref{thm:Ore1960}, 
guarantees the existence of a partition of a graph into a prescribed number of cycles.

\begin{Thm}[Brandt et al. \cite{BCFGL1997}]
\label{thm:BCFGL1997}
Let $k$ be a positive integer, 
and let $G$ be a graph of order $n \ge 4k-1$. 
If $\sigma_{2} \ge n$, 
then $G$ can be partitioned into $k$ cycles, 
i.e., 
$G$ contains $k$ disjoint cycles $C_{1}, \dots, C_{k}$ 
such that $\bigcup_{1 \le p \le k}V(C_{p}) = V(G)$. 
\end{Thm}

In order to generalize results on Hamilton cycles, 
degree conditions 
for partitioning graphs into a prescribed number of cycles with some additional conditions, 
have been extensively studied. 
See a survey paper \cite{CY2018}.

On the other hand, 
Hajnal and Szemer\'{e}di (1970) gave the following minimum degree condition 
for graphs to be partitioned into $k$ complete graphs of order $t$.

\begin{Thm}[Hajnal and Szemer\'{e}di \cite{HS1970}]
\label{thm:HS1970}
Let $k$ and $t$ be integers with $k \ge 1$ and $t \ge 3$, 
and let $G$ be a graph of order $n = tk$. 
If 
$\delta \ge \frac{t-1}{t}n$, 
then $G$ can be partitioned into $k$ subgraphs, 
each of which is isomorphic to $K_{t}$. 
\end{Thm}

In 2008, Kierstead and Kostochka improved 
the $\delta$ condition into the following $\sigma_{2}$ condition.

\begin{Thm}[Kierstead and Kostochka \cite{KK2008}]
\label{thm:KK2008}
Let $k$ and $t$ be integers with $k \ge 1$ and $t \ge 3$, 
and let $G$ be a graph of order $n = tk$. 
If 
$\sigma_{2} \ge \frac{2(t-1)}{t}n - 1$, 
then $G$ can be partitioned into $k$ subgraphs, 
each of which is isomorphic to $K_{t}$. 
\end{Thm}

The above two theorems concern with the existence of 
an equitable (vertex-)coloring in graphs. 
In fact, Theorem~\ref{thm:HS1970} 
implies that 
a conjecture of Erd\H{o}s~\cite{Erdos1964} 
(``every graph of maximum degree at most $k-1$ has an equitable $k$-coloring'') 
is true. 
Motivated by this conjecture, 
Seymour~\cite{Seymour1974} 
proposed a more general conjecture, 
which states that 
every graph of order $n \ge 3$ 
and of minimum degree at least $\frac{t-1}{t}n$ 
contains $(t-1)$-th power of a Hamilton cycle. 
It is also a generalization of Theorem~\ref{thm:Dirac1952} 
by including the case $t = 2$. 
In \cite{KSS1998}, Koml\'os et al.~proved the Seymour's conjecture for sufficiently large graphs by using the Regularity Lemma. 
For other related results, see a survey paper \cite{KKY2009}.

In this paper, we focus on a relaxed structure of a complete subgraph in graphs as follows. 
For an integer $c \ge 1$, 
a cycle $C$ in a graph $G$ is called a $c$-\textit{chorded cycle} 
if there are at least $c$ edges between the vertices on the cycle $C$ that are not edges of $C$, 
i.e., $\big| E(G[V(C)]) \setminus E(C) \big| \ge c$, 
where 
for a vertex subset $X$ of $G$, 
$G[X]$ denotes the subgraph of $G$ induced by $X$. 
We call each edge of $E(G[V(C)]) \setminus E(C)$ a \textit{chord} of $C$. 
Since a Hamilton cycle of $K_{t}$ 
has exactly $\frac{t(t-3)}{2}$ chords, 
we can regard a $c$-chorded cycle 
as a relaxed structure of $K_{t}$ for $c = \frac{t(t-3)}{2}$. 
Concerning the existence of a partition into such structures, 
we give the following result. 
Here, for positive integers $k$ and $c$, we define $f(k, c) = 8k^{2}c + 10kc - 4k + 2c + 1$.

\begin{thm}
\label{thm:degree condition for partitions into c-chorded cycles}
Let $k$ and $c$ be positive integers, 
and let $G$ be a graph of order $n \ge f(k, c)$. 
If 
$\sigma_{2} \ge n$, 
then $G$ can be partitioned into $k$ $c$-chorded cycles. 
\end{thm}

This theorem says that for a sufficiently large graph, 
the Ore condition 
also guarantees the existence of a partition into $k$ subgraphs, 
each of which is a relaxed structure of a complete graph. 
The complete bipartite graph $K_{\frac{n-1}{2}, \frac{n+1}{2}}$ ($n$ is odd) 
shows the sharpness of the lower bound on the degree condition. 
But 
we do not know whether the order condition (the function $f(k, c)$) is sharp or not. 
It comes from our proof techniques.

Related results can be found in \cite{BD2009, BLR2018, CGHOS2015, GHH2015, GHM2014}. 
In these papers, 
degree conditions for \textit{packing cycles with many chords} in a graph, 
i.e., finding a prescribed number of disjoint cycles with many chords
(it may not form a partition of a graph), 
are given 
and some of the results are also generalizations of Theorem~\ref{thm:HS1970}.

In Section~\ref{Lemmas}, 
we give lemmas which are obtained from arguments for hamiltonian problems. 
By using the lemmas, in Section~\ref{Proof of main}, 
we first show that 
the collection of disjoint $c$-chorded cycles 
in a graph $G$ satisfying the conditions of Theorem~\ref{thm:degree condition for partitions into c-chorded cycles}, 
can be transformed into a partition of $G$ (Theorem~\ref{thm:packing to partition} in Section~\ref{Proof of main}). 
Then we show that 
Theorem~\ref{thm:packing to partition} 
and a result on packing cycles 
lead to 
Theorem~\ref{thm:degree condition for partitions into c-chorded cycles} as a corollary 
(see the last of Section~\ref{Proof of main}). 
In Section~\ref{Concluding remarks}, 
we give some remarks on the order condition 
and show that 
the order condition in Theorem~\ref{thm:degree condition for partitions into c-chorded cycles} 
can be improved for the case of the Dirac condition.

\section{Lemmas}
\label{Lemmas}

We prepare terminology and notations which will be used in our proofs. 
Let $G$ be a graph. 
For $v \in V(G)$ and $X \subseteq V(G)$, 
we let $N_{X} (v) = N_{G}(v) \cap X$ 
and $d_{X} (v) = | N_{X} (v) |$. 
For $V, X \subseteq V(G)$, 
let $N_{X} (V) = \bigcup_{v \in V}N_{X} (v)$. 
For a subgraph $F$ of $G$, 
we define $\ol{E_{G}}(F) = E(G[V(F)]) \setminus E(F)$. 
A $(u, v)$-path in $G$ is a path from a vertex $u$ to a vertex $v$ in $G$. 
We write a cycle (or a path) $C$ with a given orientation by $\ora{C}$. 
If there exists no fear of confusion, 
we abbreviate $\ora{C}$ by $C$. 
Let $C$ be an oriented cycle (or path). 
We denote by $\ola{C}$ the cycle (or the path) $C$ with the reverse orientation. 
For $v \in V(C)$, 
we denote by $v^{+}$ and $v^{-}$ the successor and the predecessor of $v$ on $\ora{C}$, 
respectively. 
For $X \subseteq V(C)$, 
we define $X^{+} = \{v^{+} : v \in X\}$ and $X^{-} = \{v^{-} : v \in X\}$.
For $u, v \in V(C)$,
we denote by $C[u, v]$
the $(u, v)$-path on $\ora{C}$.
The reverse sequence of $C[u, v]$ is denoted by $\ola{C}[v, u]$. 
In the rest of this paper, 
we often identify a subgraph $F$ of $G$ with its vertex set $V(F)$.

We next prepare some lemmas. 
In the proof, we use the technique for proofs concerning hamiltonian properties of graphs.
To do that, 
in the rest of this section, 
we fix the following. 
Let $k$ and $c$ be positive integers, 
and let $G$ be a graph of order $n$ 
and $L$ a fixed vertex subset of $G$. 
Let $C_{1}, \dots, C_{k}$ be $k$ disjoint $c$-chorded cycles each with a fixed orientation in $G$, 
and suppose that $C^{*} := \bigcup_{1 \le p \le k}C_{p}$ is not a spanning subgraph of $G$. 
Let $H^{*} = G - C^{*}$ 
and $H$ be a component of $H^{*}$. 
Assume that $C_{1}, \dots, C_{k}$ are chosen so that 
\begin{enumerate}
\renewcommand{\labelenumi}{\upshape{(A\arabic{enumi}})}
\setcounter{enumi}{0}
\item 
$|V(C^{*}) \cap L|$ is as large as possible, and 

\item 
$|C^{*}|$ is as large as possible, subject to (A1). 
\end{enumerate}
Then the choices lead to the following.

\begin{lem}
\label{lem:crossing}
Let $C = C_{p}$ with $1 \le p \le k$, and let $v \in N_{C}(H)$ and $x \in V(H)$. 
Then~(i)~$v^{+}x \notin E(G)$, 
and 
(ii)~$d_{H^{*} \cup C}(v^{+}) + d_{H^{*} \cup C}(x) \le |H^{*} \cup C| - 1$. 
\end{lem}
\noindent
\textbf{Proof of Lemma~\ref{lem:crossing}.}~We let $\ora{P}$ be a $(v^{+}, x)$-path 
consisting of 
the path $C[v^{+}, v]$ 
and a $(v, x)$-path in $G[V(H) \cup \{v\}]$. 

\begin{figure}[H]
\begin{center}
\includegraphics[scale=1.00,clip]{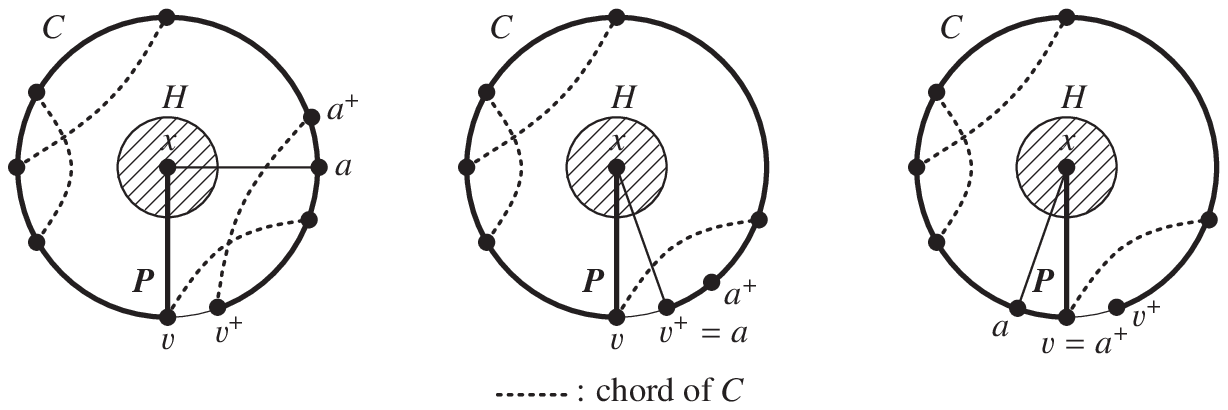}
\caption{Lemma~\ref{lem:crossing}}
\label{crossingfig}
\end{center}
\end{figure}

Suppose first that 
there exists a vertex $a$ in $\big( N_{P}(v^{+}) \big)^{-} \cap N_{P}(x)$, 
where the superscript~$^{-}$~refers to the orientation of $\ora{P}$ 
(see Figure~\ref{crossingfig}). 
Consider the cycle $D := v^{+} P[a^{+}, x] \ola{P}[a, v^{+}]$. 
Then by the definitions of $P$ and $D$, 
we have 
$\big( \ol{E_{G}}(C) \setminus \{ v^{+}a^{+} \} \big) \cup \{vv^{+}\} \subseteq \ol{E_{G}}(D)$ 
or 
$\ol{E_{G}}(C) \cup \{aa^{+}\} \subseteq \ol{E_{G}}(D)$, 
and hence $D$ is a $c$-chorded cycle in $G[V(H^{*} \cup C)]$. 
Moreover we also have 
$V(C) \subset V(P) = V(D)$. 
Therefore, 
by replacing $C$ with $D$, 
this contradicts (A1) or (A2). 
Thus 
\begin{align}
\label{crossing}
\big( N_{P}(v^{+}) \big)^{-} \cap N_{P}(x) = \emptyset. 
\end{align}
This in particular implies that 
$v^{+}x \notin E(G)$. 
Thus (i) holds. 

Suppose next that 
there exists a vertex $b$ in 
$N_{G}(v^{+}) \cap N_{G}(x) \cap \big( V(H^{*} \cup C) \setminus V(P) \big)$. 
Consider the cycle $D' := P[v^{+}, x] b v^{+}$. 
Then by the similar argument as above, 
replacing $C$ with $D'$ would violate (A1) or (A2).
Thus $N_{G}(v^{+}) \cap N_{G}(x) \cap \big( V(H^{*} \cup C) \setminus V(P) \big) = \emptyset$. 
Combining this with (\ref{crossing}), 
we get $d_{H^{*} \cup C}(v^{+}) + d_{H^{*} \cup C}(x) \le |H^{*} \cup C| - 1$. 
Thus (ii) holds. 
\qed

\vspace{-18pt}

\begin{lem}
\label{lem:crossing2}
Let $C = C_{p}$ with $1 \le p \le k$, and let $u^{*}, v^{*} \in \big( N_{C}(H) \big)^{-} \cup \big( N_{C}(H) \big)^{+}$ and $x \in V(H)$. 
If $\sigma_{2}(G) \ge n$, 
then the following hold. 
\begin{enumerate}[{\upshape(i)}]
\item $d_{C_{q}}(u^{*}) + d_{C_{q}}(x) \ge |C_{q}| + 1$ for some $q$ with $q \neq p$. 
\item $d_{C_{q'}}(u^{*}) + d_{C_{q'}}(v^{*}) \ge |C_{q'}| + 1$ for some $q'$ with $q' \neq p$. 
\end{enumerate}
\end{lem}
\noindent
\textbf{Proof of Lemma~\ref{lem:crossing2}.}~Note that by Lemma~\ref{lem:crossing}-(i)
\footnote{We use the symmetry of $\ora{C}$ and $\ola{C}$ for a vertex of $\big( N_{C}(H) \big)^{-}$.}, 
$u^{*}x, v^{*}x \notin E(G)$. 
Since $\sigma_{2}(G) \ge n$, 
it follows from Lemma~\ref{lem:crossing}-(ii)
\footnotemark[1] 
that 
\begin{align}
\label{degree sum 1}
&d_{C^{*} - C}(u^{*}) + d_{C^{*} - C}(x) \ge n - \big( | H^{*} \cup C | - 1\big) = | C^{*} - C | + 1, \\
\label{degree sum 2}
&d_{C^{*} - C}(v^{*}) + d_{C^{*} - C}(x) \ge n - \big( | H^{*} \cup C | - 1\big) = | C^{*} - C | + 1. 
\end{align}
Then (\ref{degree sum 1}) and the Pigeonhole Principle yield that (i) holds. 
Since 
$d_{C_{r}}(x) \le |C_{r}|/2$ for $1 \le r \le k$ by Lemma~\ref{lem:crossing}-(i), 
combining (\ref{degree sum 1}) and (\ref{degree sum 2}), 
and the Pigeonhole Principle yield that (ii) holds. 
\qed

\section{Proof of Theorem~\ref{thm:degree condition for partitions into c-chorded cycles}}
\label{Proof of main}

In order to show Theorem~\ref{thm:degree condition for partitions into c-chorded cycles}, 
we first prove the following theorem.

\begin{thm}
\label{thm:packing to partition}
Let $k, c$ and $G$ be the same as the ones in Theorem~\ref{thm:degree condition for partitions into c-chorded cycles}. 
Suppose that 
$G$ contains $k$ disjoint $c$-chorded cycles. 
If 
$\sigma_{2} \ge n$, 
then $G$ can be partitioned into $k$ $c$-chorded cycles. 
\end{thm}

In the proof of Theorem~\ref{thm:packing to partition}, 
we use the following lemma.

\begin{Lem}[see Lemma~2.3 in \cite{Egawa1989}]
\label{lem:Egawa1989}
Let $d$ be an integer, and 
let $G$ be a 2-connected graph of order $n$ and $a \in V(G)$. 
If $d_{G}(u) + d_{G}(v) \ge d$ for any two distinct non-adjacent vertices $u, v$ of $V(G) \setminus \{a\}$, 
then $G$ contains a cycle of order at least $\min \{d, n\}$. 
\end{Lem}

\noindent
\textbf{Proof of Theorem~\ref{thm:packing to partition}.}~Let 
$L$, $C_{1}, \dots, C_{k}$, $C^{*}$ and $H^{*}$ 
be the same ones as in the paragraph preceding Lemma~\ref{lem:crossing} in Section~\ref{Lemmas}.

\begin{claim}
\label{claim:|N_{G}(H; C)| <= 2c}
If $H$ is a component of $H^{*}$, then 
$|N_{C_{p}} (H)| \le 2c$ for $1 \le p \le k$. 
\end{claim}
\proof 
Let $H$ be a component of $H^{*}$. 
It suffices to consider the case $p = 1$. 
Suppose that $|N_{C_{1}} (H)| \ge 2c+1$. 
Let $e_{1}, \dots, e_{c}$ be $c$ distinct edges 
in $\ol{E_{G}}(C_{1})$. 
Note that by Lemma~\ref{lem:crossing}-(i), 
$N_{C_{1}} (H) \cap \big( N_{C_{1}} (H) \big)^{+} = \emptyset$. 
Then, 
since $|N_{C_{1}} (H)| \ge 2c+1$, 
we can take two distinct vertices $v_{1}, v_{2}$ in $N_{C_{1}} (H)$ 
such that 
\begin{align}
\label{chord of C_{1}}
\textup{the end vertices of $e_{1}, \dots, e_{c}$ do not appear in $C_{1}[ v_{1}^{+}, v_{2}^{-} ]$, 
i.e., 
$e_{1}, \dots, e_{c} \in \ol{E_{G}}(C_{1}[v_{2}, v_{1}])$.} 
\end{align}

We apply 
Lemma~\ref{lem:crossing2}-(ii) with $(p, u^{*}, v^{*}) = (1, v_{1}^{+}, v_{2}^{-})$. 
Then 
there exists another cycle $C_{q}$ with $q \neq 1$, say $q = 2$, 
such that 
$d_{C_{2}}( v_{1}^{+} ) + d_{C_{2}}( v_{2}^{-} ) \ge |C_{2}| + 1$. 
This inequality implies 
that 
\footnote{Change the orientation of $C_{2}$ if necessary.} 
there exists an edge $w^{-} w$ in $E(\ora{C_{2}})$ such that 
$v_{1}^{+}w^{-}, v_{2}^{-}w \in E(G)$. 
We consider two cycles 
\begin{align*}
D_{1} := C_{1}[v_{2}, v_{1}] P[v_{1}, v_{2}] 
\textup{ and }
D_{2} := C_{1}[v_{1}^{+}, v_{2}^{-}] C_{2}[w, w^{-}] v_{1}^{+},  
\end{align*}
where $P[v_{1}, v_{2}]$ denotes a $(v_{1}, v_{2})$-path in $G[V(H) \cup \{v_{1}, v_{2}\}]$ 
such that $V(P) \cap V(H) \neq \emptyset$. 
Then 
by (\ref{chord of C_{1}}), 
$D_{1}$ is a $c$-chorded cycle. 
Since 
$\ol{E_{G}}(C_{2}) \subseteq \ol{E_{G}}(D_{2})$, 
$D_{2}$ is also a $c$-chorded cycle. 
Moreover, $V(D_{1}) \cap V(D_{2}) = \emptyset$ 
and 
$V(D_{1}) \cup V(D_{2}) = V(C_{1}) \cup V(C_{2}) \cup V(P)$.  
Hence, 
replacing $C_{1}$ and $C_{2}$ with $D_{1}$ and $D_{2}$, 
this contradicts (A1) or (A2). 
\qed

Now 
we define the fixed vertex subset $L$ of $G$ as follows: 
\begin{align*}
L = \Big\{ v \in V(G) : d_{G}(v) < \frac{n}{2} \Big\}.
\end{align*}

\begin{enumerate}[{\textup{{\bf Case~\arabic{enumi}.}}}]
\setcounter{enumi}{0}
\item 
$|H^{*}| \ge \frac{n}{2} - 2kc +1$. 
\end{enumerate}

Since 
$G$ is connected, 
there exists a vertex $x \in V(H^{*})$ and a cycle $C_{p}$, say $p = 1$,  
such that $N_{C_{1}}(x) \neq \emptyset$. 
Let 
$H^{**} = H^{*} - \{x\}$.

In this case, 
we show that the following claim holds.

\begin{claim}
\label{claim:c-chorded cycle in H^{**}}
$H^{**}$ contains a $c$-chorded cycle.
\end{claim}
\proof 
We first define the following real number $\omega(c)$. 
Let $\omega(c)$ be the positive root of 
the equation $\frac{t(t-3)}{2} - c = 0$, 
i.e., $\omega(c)  = \frac{\sqrt{8c + 9} + 3}{2}$. 
Since $|E(K_{t})| - t = \frac{t(t-3)}{2}$, 
it follows that a Hamilton cycle of a complete graph of 
order at least $\lceil \omega(c) \rceil$ has at least $c$ chords.

If $V(H^{**}) \subseteq L$, then 
by the definition of $L$, 
$H^{**}$ is a complete graph, 
and hence a Hamilton cycle of $H^{**}$ has at least $c$ chords 
since $|H^{**}| \ge \frac{n}{2} - 2kc \ge \omega(c)$. 
Thus 
we may assume that $V(H^{**}) \setminus L \neq \emptyset$. 
Let $H'$ be a component of $H^{**}$ such that $V(H') \setminus L \neq \emptyset$. 
Note that by Claim~\ref{claim:|N_{G}(H; C)| <= 2c}, 
for $x' \in V(H') \setminus L$, 
$|H'| \ge d_{H'}(x') + 1 \ge \left( \frac{n}{2} - d_{C^{*}}(x') - |\{x\}| \right)+ 1 \ge \frac{n}{2} - 2kc \ge 3$.

We define an induced subgraph $B$ of $H'$ as follows: 
If $H'$ is not $2$-connected, 
let $B$ be an end block with a single cut vertex $a$ 
such that $V(B) \setminus (\{a\} \cup L) \neq \emptyset$ 
(note that we can take such a block $B$ 
because $|H'| \ge 3$ and hence $H'$ has at least two end blocks); 
If $H'$ is 2-connected, then 
let $B = H'$ and $a$ be a vertex of $H'$ such that $V(B) \setminus (\{a\} \cup L) \neq \emptyset$ 
(recall that $V(H') \setminus L \neq \emptyset$). 
Then by Claim~\ref{claim:|N_{G}(H; C)| <= 2c}, the definitions of $B$ and $a$, 
it follows that for $b \in V(B) \setminus (\{a\} \cup L)$, 
\begin{align*}
|B| \ge d_{B}(b) + 1 \ge \left( \frac{n}{2} - d_{C^{*}}(b) - |\{x\}| \right)+ 1 \ge \frac{n}{2} - 2kc. 
\end{align*}
In particular, $B$ is $2$-connected since $\frac{n}{2} - 2kc \ge 3$. 
Moreover, 
we also see that 
\begin{align*} 
\textup{$d_{B}(u) + d_{B}(v) \ge n - 4kc - 2 |\{x\}| = n - 4kc - 2$ 
for $u, v \in V(B) \setminus \{a\}$ with $u \neq v$ and $uv \notin E(G)$.} 
\end{align*}
Hence, 
by Lemma~\ref{lem:Egawa1989}, 
\begin{align*}
\textup{$B$ contains a cycle $C$ of order at least $\min \{n - 4kc - 2, |B|\}$.} 
\end{align*}
To complete the proof of the claim, 
we show that 
the cycle $C$ is a $c$-chorded cycle.

Suppose 
$G[V(C) \setminus \{a\}]$ is complete. 
Since $n \ge f(k, c)$ and $|B| \ge \frac{n}{2} - 2kc$, 
we have $|V(C) \setminus \{a\}| \ge \min \{n - 4kc - 3, |B|-1\} \ge \omega(c)$, 
and hence it follows that $C$ has at least $c$ chords. 
Thus we may assume that 
there exist two distinct non-adjacent vertices $u, v$ of $V(C) \setminus \{a\}$. 
Then by Claim~\ref{claim:|N_{G}(H; C)| <= 2c}, the definitions of $B$ and $a$, 
we have 
\begin{align*}
d_{C}(u) + d_{C}(v) 
&\ge 
n -  \big( d_{C^{*}}(u) + d_{C^{*}}(v) \big) 
-  \big( d_{B-C}(u) + d_{B - C}(v) \big) - 2|\{x\}| \\
&\ge 
n - 4kc - 2\big( |B| - |C| \big) -2\\
&\ge
n - 4kc - 2 \Big( |B| - \min \{n - 4kc - 2, |B|\} \Big) - 2\\
&= 
n - 4kc - 2 + 2 \cdot \min \{n - 4kc - 2 - |B|, 0\}. 
\end{align*}
Note that 
each $C_{i}$ has order at least $\omega(c)$ because $C_{i}$ has at least $c$ chords, 
and hence 
\begin{align*}
|B| \le |H^{**}| = |H^{*} - \{x\}| = n - 1 - |C^{*}| \le n - 1 - k \cdot \omega(c). 
\end{align*}
Since $n \ge f(k, c)$, it follows that 
\begin{align*}
d_{C}(u) + d_{C}(v) 
&\ge 
n - 4kc - 2 + 2 \cdot \min \{n - 4kc - 2 - ( n - 1 - k \cdot \omega(c) ), 0\} \\
&= 
n - 12kc + 2k \cdot \omega(c) - 4 \ge c + 4.  
\end{align*}
This implies that 
$C$ has at least $c$ chords. 
\qed

Now let 
$D_{1}$ be a $c$-chorded cycle in $H^{**} \ (= H^{*} - \{x\})$. 
Recall 
that $N_{C_{1}}(x) \neq \emptyset$. 
Let $v \in N_{C_{1}}(x)$. 
Then 
by Lemma~\ref{lem:crossing2}-(i), 
there exists a cycle $C_{q}$ with $q \neq 1$, say $q = 2$, 
such that 
$d_{C_{2}}(v^{+}) + d_{C_{2}}(x) \ge |C_{2}| + 1$. 
This inequality implies that 
there exists an edge $w^{-}w$ in $E(\ora{C_{2}})$ such that 
$v^{+}w^{-}, xw \in E(G)$. 
Let $D_{2} = C_{1}[v^{+}, v] \ x \ C_{2}[w, w^{-}] v^{+}$. 
Then, since 
$\ol{E_{G}}(C_{2}) \subseteq \ol{E_{G}}(D_{2})$, 
$D_{2}$ is a $c$-chorded cycle. 
Moreover, 
$V(D_{1}) \cap V(D_{2}) = \emptyset$ 
and 
$V(C_{1}) \cup V(C_{2}) \subset V(D_{1}) \cup V(D_{2}) \subseteq V(C_{1}) \cup V(C_{2}) \cup V(H^{*})$. 
Hence, 
replacing $C_{1}$ and $C_{2}$ with $D_{1}$ and $D_{2}$ 
would violate (A1) or (A2), a contradiction.

This completes the proof of Case~1.

\begin{enumerate}[{\textup{{\bf Case~\arabic{enumi}.}}}]
\setcounter{enumi}{1}
\item 
$|H^{*}| < \frac{n}{2} - 2kc +1$. 
\end{enumerate}

The following two claims are essential parts in this case.

\begin{claim}
\label{claim:all vertices of H are L}
\begin{enumerate}[{\upshape(i)}]
\item $V(H^{*}) \subseteq L$ (in particular, $H^{*}$ is complete) and $|H^{*}| \le 2c+1$. 
\item $\big( V(C^{*}) \setminus N_{C^{*}}(H^{*}) \big) \cap L = \emptyset$. 
\item $d_{C_{p}}(v) \ge |C_{p}| - 2kc + 1$ for $1 \le p \le k$ and $v \in V(C_{p}) \setminus N_{C_{p}}(H^{*})$. 
\end{enumerate}
\end{claim}
\proof 
We first show (i) and (ii). 
If there exists a vertex $x$ of $V(H^{*})$ such that $x \notin L$, 
then by Claim~\ref{claim:|N_{G}(H; C)| <= 2c}, 
$|H^{*}| \ge d_{H^{*}}(x) + |\{x\}| \ge \big( \frac{n}{2} - 2kc \big) + 1$, which contradicts the assumption of Case~2. 
Thus $$V(H^{*}) \subseteq L.$$ 
In particular, $H^{*}$ is a complete graph.  
Then by the definition of $L$, 
we have $$\big( V(C^{*}) \setminus N_{C^{*}}(H^{*}) \big) \cap L = \emptyset.$$ 
This together with Claim~\ref{claim:|N_{G}(H; C)| <= 2c} implies that 
$|V(C_{p}) \cap L| \le 2c$ for $1 \le p \le k$. 
Therefore, 
if $|H^{*}| \ge 2c + 2$, 
then 
by replacing the cycle $C_{1}$ with a Hamilton cycle of $H^{*}$, 
this contradicts (A1). 
Thus we have
\footnote{This argument actually implies that $|H^{*}| \le \max \{ 2c, 3 \}$. 
But we make no attempt to optimize the upper bound on $|H^{*}|$ 
since it does not lead to a significant improvement of the condition on $n$.} 
$$|H^{*}| \le 2c+1.$$

We finally show (iii). 
Let $1 \le p \le k$ and $v \in V(C_{p}) \setminus N_{C_{p}}(H^{*})$. 
We may assume that $p = 1$. 
Let $x$ be an arbitrary vertex of $H^{*}$. 
Then by Claim~\ref{claim:|N_{G}(H; C)| <= 2c}, and since $v \notin N_{C_{p}}(H^{*})$, 
we get 
\begin{align*}
d_{C^{*}}(v) \ge n - d_{C^{*}}(x) - d_{H^{*}}(x) \ge n - 2kc - (|H^{*}| - 1)= |C^{*}| - 2kc + 1. 
\end{align*}
Since 
$d_{C_{q}}(v) \le |C_{q}|$ for $2 \le q \le k$, 
we have $d_{C_{1}}(v) \ge |C_{1}| - 2kc + 1$. 
Thus (iii) holds. 
\qed

\begin{claim}
\label{claim:two chords}
Let $C = C_{p}$ with $1 \le p \le k$, 
and $w^{-} w \in E(\ora{C})$ and $S = N_{C}(H^{*})$. 
If $|C| \ge 8kc + 10c - 4$, then there exist two distinct chords $u_{1}v_{1}, u_{2}v_{2}$ of $C$ 
satisfying the following conditions (A)--(C).  
\begin{enumerate}[{\upshape(A)}]
\item $u_{1}, u_{2}, v_{2}, v_{1}$ are appear in the order along $\ora{C}$, 
\item $w^{-}, w \in C[v_{1}, u_{1}]$ and $S \subseteq C[v_{1}, u_{1}] \cup C[u_{2}, v_{2}]$, 
\item $d_{C[v_{1}, u_{1}]}(u_{1}) \ge c + 2$ and  $d_{C[u_{2}, v_{2}]}(u_{2}) \ge c + 2$. 
\end{enumerate}
\end{claim}
\proof 
Note that by Claim~\ref{claim:all vertices of H are L}-(i), 
$H^{*}$ consists of exactly one component, 
and hence 
Claim~\ref{claim:|N_{G}(H; C)| <= 2c} 
yields that $|S| \le 2c$. 
Note also that 
by Claim~\ref{claim:all vertices of H are L}-(iii), 
$d_{C}(v) \ge |C| - 2kc + 1$ for $v \in V(C) \setminus S$.

\begin{figure}[H]
\begin{center}
\hspace{-30pt}\includegraphics[scale=1.00,clip]{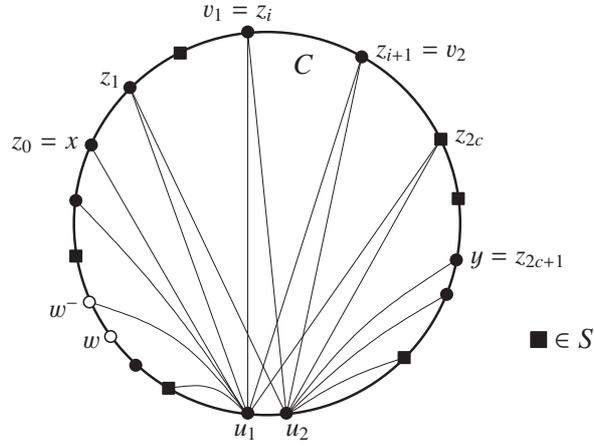}
\caption{The vertices $u_{1}, u_{2}, v_{1}, v_{2}$}
\label{twochordsfig}
\end{center}
\end{figure}

We first 
define four vertices $u_{1}, u_{2}, x, y$ of $V(C)$ 
by the following procedure (I)--(III) 
(the vertices $u_{1}, u_{2}$ will be the end vertices of the desired chords, 
and the vertices $x, y$ will be candidates of the end vertices of the desired chords). 
See also Figure~\ref{twochordsfig}. 

\begin{enumerate}[{\upshape(I)}]
\item 
Let 
$u_{1}, u_{2}$ be vertices of $V(C)$ such that 
\begin{center}
\begin{minipage}{0.45\textwidth}
\begin{equation}\label{u_{1}, u_{2}}
u_{1} = u_{2}^{-},  \tag{I-1}
\end{equation}
\end{minipage}
\hfill
\begin{minipage}{0.45\textwidth}
\begin{equation}\label{u_{1} not w, w^{-}, S}
\textup{and } \ \ u_{1}, u_{2} \notin \{w^{-}\} \cup S. \tag{I-2}
\end{equation}
\end{minipage}
\end{center}
Note that we can take such two vertices 
because $|C| \ge 8kc + 10c - 4$ and $|\{w^{-}\} \cup S| \le 2c+1$. 
Choose $u_{1}, u_{2}$ so that 
$| C[w, u_{1}] |$ is as small as possible. 
Then by the choice, 
\begin{align}
\label{from w to u_{1}}
| C[w, u_{1}] | \le 2|S|+|\{u_{1}\}|  \le 4c + 1 \tag{I-3}.
\end{align} 

\item 
Since 
$d_{C}(u_{1}) \ge |C| - 2kc + 1 \ge c + 2$ 
and $u_{1}u_{2} \in E(G)$, and by (\ref{u_{1}, u_{2}}), (\ref{u_{1} not w, w^{-}, S}), 
we can take a vertex $x$ of $N_{C}(u_{1})$ 
such that 
\begin{center}
\begin{minipage}{0.45\textwidth}
\begin{equation}
\label{x1}
w^{-} \in C[x, u_{1}], \tag{II-1}\end{equation}
\end{minipage}
\hfill
\begin{minipage}{0.45\textwidth}
\begin{equation}
\label{x2}
\textup{and } \ \ d_{C[x, u_{1}]}(u_{1}) \ge c + 2. \tag{II-2}
\end{equation}
\end{minipage}
\end{center}
In fact, the vertex $u_{2}$ can be such a vertex $x$. 
Choose $x$ so that $d_{C[x, u_{1}]}(u_{1})$ is as small as possible, 
subject to (\ref{x1}) and (\ref{x2}). 
Then by the choice, 
\begin{align*}
&\textup{if 
$d_{C[w, u_{1}]}(u_{1}) \le c+1$,}\\
&\hspace{+24pt}\textup{then $d_{C[x, w^{-}]}(u_{1}) = c + 2 - d_{C[w, u_{1}]}(u_{1})$, 
that is, 
$d_{C[x, u_{1}]}(u_{1}) = c + 2$;}\\
&\textup{if 
$d_{C[w, u_{1}]}(u_{1}) \ge c + 2$,}\\
&\hspace{+24pt}\textup{then $d_{C[x, w^{-}]}(u_{1}) = |\{x\}| = 1$, 
that is, 
$d_{C[x, u_{1}]}(u_{1}) \le |V(C[w, u_{1}]) \setminus \{u_{1}\}| + 1$.} 
\end{align*}
In either case, by (\ref{from w to u_{1}}), 
\begin{align}
\label{choice of x (1)}
d_{C[x, u_{1}]}(u_{1}) \le 4c + 1. \tag{II-3}
\end{align} 
Moreover, 
since 
$|V(C) \setminus N_{C}(u_{1}) | \le |C|  - (|C| - 2kc + 1) = 2kc - 1$, 
we have 
\begin{align}
\label{choice of x (2)}
|C[x, u_{1}]| 
&= 
\big| N_{C[x, u_{1}]}(u_{1}) \big| + \big| V(C[x, u_{1}]) \setminus N_{C[x, u_{1}]}(u_{1}) \big| \notag\\ 
&\le 
(4c + 1) + (2kc - 1) = 2kc + 4c. \tag{II-4}
\end{align}

\item 
Let $y$ be the vertex of $N_{C}(u_{2})$ 
such that 
\begin{align}
\label{y}
d_{C[u_{2}, y]}(u_{2}) = c + 2. \tag{III-1}
\end{align}
By the similar argument as in (\ref{choice of x (2)}), 
we have 
\begin{align}
\label{choice of y}
|C[u_{2}, y]| 
\le (c + 2)  + (2kc - 1) = 2kc + c + 1. \tag{III-2}
\end{align}
Recall that $|C| \ge 8kc + 10c - 4$. 
Hence by the definitions of $x, y$ and, (\ref{u_{1}, u_{2}}), (\ref{choice of x (2)}) and (\ref{choice of y}), 
\begin{align}
\label{order of y, x}
\textup{$y$ and $x$ appear in the order along $C[u_{2}^{+}, u_{1}^{-}]$, and $y^{+} \neq x$.} \tag{III-3}
\end{align}
\end{enumerate}

\noindent
To complete the proof of the claim, 
we next define two vertices $v_{1}, v_{2}$ of $V(C)$ as follows. 

\begin{enumerate}[{\upshape(IV)}]
\item 
We first show that 
\begin{align}
\label{common neighbor}
\big| N_{C[y^{+}, x^{-}]} (u_{1}) \cap N_{C[y^{+}, x^{-}]} (u_{2}) \big| \ge 2c. \tag{IV-1}
\end{align}
Assume not. 
Then 
for some $i$ with $i \in \{1, 2\}$, 
\begin{align*}
d_{C[y^{+}, x^{-}]}(u_{i}) 
&\le 
\frac{1}{2} \left( \left| V( C[y^{+}, x^{-}] ) 
\setminus 
\left( N_{C[y^{+}, x^{-}]} (u_{1}) \cap N_{C[y^{+}, x^{-}]} (u_{2}) \right) \right| \right)\\
&\hspace{+12pt}+ 
\big| N_{C[y^{+}, x^{-}]} (u_{1}) \cap N_{C[y^{+}, x^{-}]} (u_{2}) \big|\\
&\le 
\frac{1}{2} \left( | C[y^{+}, x^{-}] | - (2c-1) \right) + (2c-1). 
\end{align*}
If this inequality holds for $i = 1$, 
then by (\ref{u_{1}, u_{2}}), (\ref{x2}), (\ref{choice of x (1)}) and (\ref{y})--(\ref{order of y, x}), 
\begin{align*}
|C| - 2kc + 1 
&\le 
d_{C}(u_{1}) \\
&\le 
d_{C[x, u_{1}]}(u_{1}) + d_{C[u_{2}, y]}(u_{1}) + d_{C[y^{+}, x^{-}]}(u_{1}) \\
&\le 
(4c + 1) + (2kc + c + 1) + \frac{1}{2} \big( | C[y^{+}, x^{-}] | - (2c-1) \big) + (2c-1)\\
&=
\frac{|C|}{2} - \frac{1}{2} \big( | C[x, u_{1}] | + | C[u_{2}, y] | \big) + 2kc + 6c + \frac{3}{2} \\
&\le
\frac{|C|}{2} - \frac{1}{2} \left( c+3 + c+ 3 \right) + 2kc + 6c + \frac{3}{2} \\
&= 
\frac{|C|}{2} + 2kc + 5c - \frac{3}{2}.
\end{align*}
This implies that $|C| \le 8kc + 10c - 5$, a contradiction. 
Similarly, for the case $i = 2$, 
it follows from (\ref{u_{1}, u_{2}}), (\ref{x2}), (\ref{choice of x (2)}), (\ref{y}) and (\ref{order of y, x}) that 
$|C| \le 8kc + 10c - 5$, a contradiction again. 
Thus (\ref{common neighbor}) is proved.

By (\ref{common neighbor}), 
we can take 
$2c$ distinct vertices $z_{1}, \dots, z_{2c}$ in $N_{C[y^{+}, x^{-}]} (u_{1}) \cap N_{C[y^{+}, x^{-}]} (u_{2})$. 
We may assume that 
$z_{2c}, \dots, z_{2}, z_{1}$ appear in the order along $C[y^{+}, x^{-}]$. 
Let $z_{0} = x$ and $z_{2c + 1} = y$. 
Then 
\begin{align}
\label{order of z}
\textup{$z_{2c + 1}, z_{2c}, \dots, z_{2}, z_{1}, z_{0}$ appear in the order along $C[u_{2}^{+}, u_{1}^{-}]$.} \tag{IV-2}
\end{align}
Note also that 
\begin{align}
\label{z_{i} is a neighbor of u_{1} and u_{2}}
\textup{$z_{i} \in N_{G}(u_{1})$ for $0 \le i \le 2c$ 
and 
$z_{i} \in N_{G}(u_{2})$ for $1 \le i \le 2c+1$.} \tag{IV-3}
\end{align}
Moreover, 
since 
$|S| \le 2c$, 
it follows that 
there exists an index $i$ with $0 \le i \le 2c$ 
such that 
\begin{align}
\label{no vertex of S}
\textup{$z_{i} = z_{i+1}^{+}$, 
or 
$z_{i} \neq z_{i+1}^{+}$ and $C[z_{i+1}^{+}, z_{i}^{-}] \cap S = \emptyset$.} \tag{IV-4}
\end{align}
Then we define
\begin{align}
\label{definition of v_{2}, v_{2}}
\textup{$v_{1} = z_{i}$ and $v_{2} = z_{i+1}$.} \tag{IV-5}
\end{align}
\end{enumerate}

Now let $u_{1}, u_{2}, v_{1}, v_{2}$ be the vertices defined as in the above (I)--(IV). 
By (\ref{z_{i} is a neighbor of u_{1} and u_{2}}) and (\ref{definition of v_{2}, v_{2}}), 
$u_{1}v_{1}$ and $u_{2}v_{2}$ are chords of $C$. 
By (\ref{order of z}) and (\ref{definition of v_{2}, v_{2}}), 
we also see that 
$u_{1}, u_{2}, v_{2}, v_{1}$ appear in the order along $\ora{C}$. Thus (A) holds. 
By (\ref{u_{1}, u_{2}}), (\ref{u_{1} not w, w^{-}, S}), (\ref{x1}), (\ref{no vertex of S}) and (\ref{definition of v_{2}, v_{2}}), 
we have 
$w^{-}, w \in C[v_{1}, u_{1}]$ and 
$S \subseteq C[v_{1}, v_{2}] = C[v_{1}, u_{1}] \cup C[u_{2}, v_{2}]$. 
Thus (B) holds. 
By (\ref{x2}), (\ref{y}) and (\ref{definition of v_{2}, v_{2}}), 
we have 
$d_{C[v_{1}, u_{1}]}(u_{1}) \ge c + 2$ and  $d_{C[u_{2}, v_{2}]}(u_{2}) \ge c + 2$. 
Thus (C) also holds. 

This completes the proof of Claim~\ref{claim:two chords}. 
\qed

Let $x \in V(H^{*})$ and $C_{p}$ be a cycle with $1 \le p \le k$ 
such that $N_{C_{p}}(x) \neq \emptyset$. 
Let $v \in N_{C_{p}}(x)$. 
We may assume that $p = 1$. 
Then 
by Lemma~\ref{lem:crossing2}-(i), 
there exists a cycle $C_{q}$ with $q \neq 1$, say $q = 2$, 
such that 
$d_{C_{2}}(v^{+}) + d_{C_{2}}(x) \ge |C_{2}| + 1$. 
This inequality implies that 
there exists an edge $w^{-}w$ in $E(\ora{C_{2}})$ such that 
$v^{+}w^{-}, xw \in E(G)$. 
On the other hand, 
since 
$|C^{*}| = n - |H^{*}| \ge n - 2c - 1$ by Claim~\ref{claim:all vertices of H are L}-(i), 
there exists a cycle $C_{r}$ with $1 \le r \le k$ 
such that $|C_{r}| \ge \frac{1}{k}(n - 2c - 1)$.

Suppose that $r \ge 3$, say $r = 3$. 
Then, since 
$|C_{3}| \ge \frac{1}{k}(n - 2c - 1) \ge \frac{1}{k}\big( f(k, c) - 2c - 1 \big) \ge 8kc + 10c - 4$, 
we can apply Claim~\ref{claim:two chords} 
to $C_{3}$ with $S = N_{C_{3}}(H^{*})$ \ \footnote{We do not use $w^{-}w$ in Claim~\ref{claim:two chords}.},  
i.e., 
$C_{3}$ has two chords $u_{1}v_{1}, u_{2}v_{2}$ satisfying the conditions (A)--(C). 
Let 
\begin{align*}
D_{1} := C_{1}[v^{+}, v] \ x \ C_{2}[w, w^{-}] v^{+}, \ 
D_{2} := u_{1} C_{3}[v_{1}, u_{1}] 
\textup{ and }  
D_{3} := C_{3}[u_{2}, v_{2}] u_{2}.  
\end{align*}
Since 
$\ol{E_{G}}(C_{1}) \subseteq \ol{E_{G}}(D_{1})$, 
$D_{1}$ is a $c$-chorded cycle. 
By the condition (C), 
$D_{2}$ and $D_{3}$ are also $c$-chorded cycles. 
By the definitions of $D_{1}, D_{2}, D_{3}$, 
the condition (B) 
and Claim~\ref{claim:all vertices of H are L}-(ii), 
we have 
$V(C_{1}) \cup V(C_{2}) \cup (V(C_{3}) \cap L) \cup \{x\} 
\subseteq 
\bigcup_{1 \le s \le 3}V(D_{s}) 
\subseteq 
\bigcup_{1 \le s \le 3}V(C_{s}) \cup \{x\}$. 
Moreover, 
by the condition (A), 
$D_{1}, D_{2}$ and $D_{3}$ are disjoint. 
Since $x \in L$ by Claim~\ref{claim:all vertices of H are L}-(i), 
replacing $C_{1}, C_{2}$ and $C_{3}$ with $D_{1}, D_{2}$ and $D_{3}$ 
would violate (A1), a contradiction.

Suppose next that $r \in \{1, 2\}$, say 
\footnote{Since $(C_{1}, v, v^{+})$ and $(\ola{C_{2}}, w, w^{-})$ are symmetric, 
we may assume that $r = 2$.}
$r = 2$. 
We apply Claim~\ref{claim:two chords} to $C_{2}$ 
so that the edge $w^{-}w$ of $C_{2}$ is the same one as in Claim~\ref{claim:two chords} 
and $S = N_{C_{2}}(H^{*})$, 
i.e., 
$C_{2}$ has two chords $u_{1}v_{1}, u_{2}v_{2}$ satisfying the conditions (A)--(C). 
Let 
\begin{align*}
D_{1} := C_{1}[v^{+}, v] \ x \ C_{2}[ w, u_{1}] C_{2}[v_{1}, w^{-}] v^{+}\textup{ and }  
D_{2} := C_{2}[u_{2}, v_{2}] u_{2}.  
\end{align*}
Since 
$\ol{E_{G}}(C_{1}) \subseteq \ol{E_{G}}(D_{1})$, 
$D_{1}$ is a $c$-chorded cycle. 
By the condition (C), 
$D_{2}$ is also a $c$-chorded cycle. 
By the definitions of $D_{1}, D_{2}$, 
the condition (B) 
and Claim~\ref{claim:all vertices of H are L}-(ii), 
we have 
$V(C_{1}) \cup (V(C_{2}) \cap L) \cup \{x\} 
\subseteq 
V(D_{1}) \cup V(D_{2}) 
\subseteq 
V(C_{1}) \cup V(C_{2}) \cup \{x\} $. 
Moreover, 
by the condition (A), 
$D_{1}$ and $D_{2}$ are disjoint. 
Since $x \in L$, 
replacing $C_{1}$ and $C_{2}$ with $D_{1}$ and $D_{2}$
would violate (A1), a contradiction again.

This completes the proof of Theorem~\ref{thm:packing to partition}. 
\qed

\medskip
We finally prove Theorem~\ref{thm:degree condition for partitions into c-chorded cycles}. 
In 2009, 
Babu and Diwan gave the following result concerning 
the existence of $k$ disjoint $c$-chorded cycles in graphs. 
(They actually proved a stronger result, see \cite{BD2009} for the detail. 
See also \cite[Theorem~3.4.16]{CY2018}.)

\begin{Thm}[Babu and Diwan \cite{BD2009}]
\label{thm:BD2009}
Let $k$ and $c$ be positive integers, 
and let $G$ be a graph of order at least $k(c+3)$. 
If 
$\sigma_{2} \ge 2k(c+2) - 1$, 
then $G$ contains $k$ disjoint $c$-chorded cycles. 
\end{Thm}

Combining this with Theorem~\ref{thm:packing to partition}, 
we get Theorem~\ref{thm:degree condition for partitions into c-chorded cycles} as follows.

\medskip
\noindent
\textbf{Proof of Theorem~\ref{thm:degree condition for partitions into c-chorded cycles}.}~Let 
$k, c$ and $G$ be the same ones as in Theorem~\ref{thm:degree condition for partitions into c-chorded cycles}, 
and suppose $\sigma_{2}(G) \ge n$. 
Since $\sigma_{2}(G) \ge n \ge f(k, c) \ge \max\{k(c+3), 2k(c+2) - 1\}$, 
Theorem~\ref{thm:BD2009} yields that 
$G$ contains $k$ disjoint $c$-chorded cycles. 
Then by Theorem~\ref{thm:packing to partition}, 
$G$ can be partitioned into $k$ $c$-chorded cycles. 
\qed

\section{Concluding remarks}
\label{Concluding remarks}

In this paper, 
we have shown that 
for a sufficiently large graph $G$, 
the Ore condition for partitioning the graph $G$ into $k$ cycles (Theorem~\ref{thm:BCFGL1997}), 
also guarantees the existence of a partition of $G$ into $k$ cycles with $c$ chords 
which are relaxed structures of a complete graph 
(see Theorem~\ref{thm:degree condition for partitions into c-chorded cycles}). 
But, as mentioned in Section~\ref{Introduction}, 
we do not know whether the order condition (the function $f(k, c)$) is sharp or not. 
Perhaps, a weaker order condition may suffice to guarantee the existence.

For the case of the Dirac condition, 
it follows from our arguments that 
the order condition can be improved as follows. 
If we assume $\delta(G) \ge \frac{n}{2}$, 
then we have $L = \emptyset$ in the proof of Theorem~\ref{thm:packing to partition}, 
i.e.,  
Case~2 does not occur (see Claim~\ref{claim:all vertices of H are L}-(i)). 
On the other hand, 
in the proof of Case~1 of Theorem~\ref{thm:packing to partition}, 
we have used the order condition in the following parts: 
\begin{itemize}
\setlength{\parskip}{0cm} 
\setlength{\itemsep}{0cm}
\item $\frac{n}{2} - 2kc \ge \omega(c) \ (= \frac{\sqrt{8c + 9} + 3}{2} \ge 3)$, 
\item $\min \{n - 4kc - 3, |B|-1\} \ge \min \{n - 4kc - 3, \frac{n}{2} - 2kc-1\} \ge \omega(c)$, 
\item $n - 12kc + 2k \cdot \omega(c) - 4 \ge c + 4$. 
\end{itemize}
In the proof of Theorem~\ref{thm:degree condition for partitions into c-chorded cycles}, 
we have also used the order condition in the following part: 
\begin{itemize}
\setlength{\parskip}{0cm} 
\setlength{\itemsep}{0cm}
\item $n \ge \max\{k(c+3), 2k(c+2) - 1\}$. 
\end{itemize}
Therefore, 
as a corollary of our arguments, 
we get the following.

\begin{thm}
\label{thm:min degree condition for partitions into c-chorded cycles}
Let $k$ and $c$ be positive integers, 
and let $G$ be a graph of order $n  \ge 12kc -  2k \cdot \omega(c) +  c + 8$, 
where $\omega(c) = \frac{\sqrt{8c + 9} + 3}{2}$. 
If 
$\delta \ge \frac{n}{2}$, 
then $G$ can be partitioned into $k$ $c$-chorded cycles. 
\end{thm}

We finally remark about the necessary order condition. 
Let $c$ be a positive integer, 
and 
let $\psi(c)$ be the positive root of 
the equation $t(t-2) - c = 0$, 
i.e., $\psi(c)  = \sqrt{c + 1} + 1$. 
Note that $|E(K_{t, t})| - 2t = t(t-2)$. 
If 
a bipartite graph contains a $c$-chorded cycle, 
then 
by the definition of $\psi(c)$, 
it follows that the order of the bipartite graph is at least $2 \lceil \psi(c) \rceil$. 
Therefore, 
the complete bipartite graph $G \cong K_{k \lceil \psi(c) \rceil -1, k \lceil \psi(c) \rceil -1}$ 
satisfies $\delta(G) = |G|/2$ and $\sigma_{2}(G) = |G|$, but $G$ cannot be partitioned into $k$ $c$-chorded cycles. 
Thus the order at least $2k \lceil \psi(c) \rceil -1$ is necessary, 
and 
the order conditions in Theorems~\ref{thm:degree condition for partitions into c-chorded cycles} 
and \ref{thm:min degree condition for partitions into c-chorded cycles} 
might be improved into $n \ge 2k \lceil \psi(c) \rceil -1$. 
Theorem~\ref{thm:BCFGL1997} supports it by including the case $c = 0$, since $\psi(c) = 2$ for the case $c = 0$. 
For the case $c = 1$, 
related results can be also found in \cite[Corollary~3.4.7]{CY2018}.



\end{document}